\documentclass[a4paper,11pt]{article}
\usepackage[active]{srcltx}
\usepackage{amsmath}\usepackage{amssymb}\usepackage{amscd}

\newcommand{\tj}{\tilde{j}}
\newcommand{\bs}[1]{{\sf{\mathbf{s}}}_{#1}}
\newcommand{\bt}{{\sf{\mathbf{t}}}}
\newcommand{\bj}[1]{{\sf{\mathbf j}}_{#1}}
\newcommand{\cc}{c}

\newcounter{pic}

\newtheorem{defi}{Definition}\newtheorem{theo}[defi]{Theorem}\newtheorem{prop}[defi]{Proposition}\newtheorem{lemm}[defi]{Lemma}

\begin{document}

$\ $

\begin{center}

\bigskip
{\Large\bf Fusion procedure for Coxeter groups of type B and complex\\[.2cm] reflection groups G(m,1,n)} 

\vspace{.6cm}

{\large {\bf O. V. Ogievetsky$^{\circ\diamond}$\footnote{On leave of absence from P. N. Lebedev Physical Institute, Leninsky Pr. 53,
117924 Moscow, Russia} and L. Poulain d'Andecy$^{\circ}$}}

\vskip 0.6cm

$\circ\ ${\large Center of Theoretical Physics\footnote{Unit\'e Mixte de Recherche (UMR 6207) du CNRS et des Universit\'es
Aix--Marseille I, Aix--Marseille II et du Sud Toulon -- Var; laboratoire
affili\'e \`a la FRUMAM (FR 2291)}, Luminy \\
13288 Marseille, France}

\vspace{.6cm}
$\diamond\ ${\large 
J.-V. Poncelet French-Russian Laboratory, UMI 2615 du CNRS, Independent University of Moscow, 11 B. Vlasievski per., 119002 Moscow, Russia
}

\end{center}

\vskip 1cm
\begin{abstract}
A complete system of primitive pairwise orthogonal idempotents for the Coxeter groups of type $B$ and, more generally, for the complex reflection groups 
$G(m,1,n)$ is constructed by a sequence of evaluations of a rational function in several variables with values in the group ring. The evaluations
correspond to the eigenvalues of the two arrays of Jucys--Murphy elements. 
\end{abstract}

\section{{\hspace{-0.55cm}.\hspace{0.55cm}}Introduction}

A. Jucys \cite{Ju} gave a construction of a complete system of pairwise orthogonal primitive idempotents of the group ring of the symmetric group;
the construction, called now {\it fusion procedure}, involves a rational function in several variables and the idempotents are obtained by taking certain limiting values of this function. We refer to \cite{Ch,Gr1,GP,JKMO,Na2,Na3,NT} for different aspects and applications of the fusion procedure for the symmetric groups. There are analogues of the fusion procedure for the Hecke algebra of type A \cite{Gr2,Na4} and the spinor extension of the symmetric 
group \cite{Na1,Jo} (see \cite{JN} for its $q$-analogue).

A version of the fusion procedure for the symmetric group was given by A. Molev in \cite{Mo}. Here the idempotents are obtained by consecutive evaluations of the rational function. An analogue of this fusion procedure was developed for the Hecke algebra \cite{IMO}, the Brauer algebra \cite{IM,IMO2} 
and the Birman--Murakami--Wenzl algebra \cite{IMO3}. 

The aim of this paper is to give a fusion procedure, in the spirit of \cite{Mo}, for the complex reflection groups $G(m,1,n)$. As in  \cite{Mo}, and later 
\cite{IM,IMO,IMO2,IMO3}, we use the Jucys--Murphy elements. They were introduced for $G(m,1,n)$ independently in \cite{Pu} and \cite{Wa}. The   
Jucys--Murphy elements of $G(m,1,n)$ form two arrays, $j_i$ and $\tilde{j}_i$, $i=1,\dots ,n$, and their union is the maximal commutative set in 
$\mathbb{C}G(m,1,n)$, see \cite{OP,Pu}. An irreducible representation of the group $G(m,1,n)$ is coded by an $m$-tuple of partitions
and the elements of the semi-normal basis correspond to standard $m$-tuples of tableaux; the eigenvalues of  $j_i$ carry information about 
the ``position" - the place of a  tableau in the $m$-tuple - while the eigenvalues of $\tilde{j}_i$ are related to the classical contents of nodes. 
In the work \cite{OP} both sets appeared as classical limits of simple expressions involving the single set of the Jucys--Murphy elements of the cyclotomic Hecke algebra, the flat deformation of $\mathbb{C}G(m,1,n)$. By the maximality, all diagonal matrix units of $\mathbb{C}G(m,1,n)$ can be expressed in terms of the Jucys--Murphy elements $j_i$ and $\tj_i$, $i=1,\dots,n$. Then we translate this expression as a fusion 
procedure:  any diagonal matrix unit can be obtained by a sequence of evaluations of
a certain rational function with values in $\mathbb{C}G(m,1,n)$. The arrays $j_i$ and $\tilde{j}_i$ play different roles: the positions can be evaluated 
simultaneously while the contents should then be evaluated consequently from 1 to $n$. 

The group $G(1,1,n)$ is isomorphic to the symmetric group $S_n$ and our fusion procedure for $m=1$ reproduces the fusion procedure of \cite{Mo}.

The group $G(2,1,n)$ is isomorphic to the hyperoctahedral group ${\sf B}_n$, the Coxeter group of type B. Thus, in particular, we obtain a fusion procedure for the Coxeter group of type $B$ and a description of a complete set of pairwise orthogonal primitive idempotents of ${\sf B}_n$ in terms 
of a single rational function with values in the group algebra $\mathbb{C}{\sf B}_n$.

For the clarity of the exposition, we first describe the fusion procedure for the Coxeter group ${\sf B}_n$. This is done not only for aesthetic reasons: 
the rational function with values in $\mathbb{C}{\sf B}_n$ leading to the complete set of idempotents can be viewed as a word for the longest 
element of ${\sf B}_n$ in which certain entries are ``Baxterized", similarly to the rational function for the type A. Although $G(m,1,n)$ is not a Coxeter group for $m>2$ and the notion of length of an element is not defined, there is an analogue of the longest
element: it is longest with respect to the normal form of \cite{OP} (the classical limit of the normal form for the cyclotomic Hecke algebra $H(m,1,n)$ \cite{OPdA1}); the rational  function with values in $\mathbb{C}G(m,1,n)$ leading to the complete set of idempotents can again be viewed as a word for the longest element with certain entries Baxterized. 
We stress that the groups $G(m,1,n)$ admit a fusion procedure for any positive integer $m$, and our construction is uniform for all $m$.

The paper is organized as follows.
Section \ref{sec-def} contains necessary definitions and notations  about the groups ${\sf B}_n$ and their representations. 
The Jucys--Murphy elements for ${\sf B}_n$ were defined in \cite{Pu,Ra,Wa}; we describe the diagonal matrix units in terms of them.
In Section \ref{sec-fus-B} we prove the Theorem \ref{prop-fus} which gives the fusion procedure for the groups ${\sf B}_n$.
In Section \ref{sec-fus-G}, in the Theorem \ref{prop-fus-m}, we extend the results 
to the general case of the groups $G(m,1,n)$. As the proofs work mainly along the same lines as for ${\sf B}_n$, we 
only indicate the necessary modifications.
 
\section{{\hspace{-0.55cm}.\hspace{0.55cm}}Idempotents and Jucys--Murphy elements of the group ${\sf B}_n$}\label{sec-def}

\subsection{Definitions}

The Coxeter group ${\sf A}_n$ of type $A$ (the symmetric group on $n+1$ letters) is generated by the elements $s_1,\dots,s_n$ with 
the defining  relations:
\begin{equation}\label{rel-def-An}
\begin{array}{ll}
s_i^2=1 &\textrm{for $i=1,\dots,n$,}\\[0.2em]
s_is_{i+1}s_i=s_{i+1}s_is_{i+1} &\textrm{for $i=1,\dots,n-1$,}\\[0.2em]
s_is_j=s_js_i &\textrm{for $i,j=1,\dots,n$ such that $|i-j|>1$.}
\end{array}
\end{equation}
The Coxeter group ${\sf B}_{n+1}$ of type $B$ (also called the hyperoctahedral group) is generated by the elements $s_1,\dots,s_n$ and  
$t$ with the defining relations (\ref{rel-def-An}), 
\begin{equation}\label{rel-def-Bn}
\begin{array}{ll}
ts_1ts_1=s_1ts_1t,\\[0.2em]
s_it=ts_i &\textrm{for $i=2,\dots,n$}
\end{array}
\end{equation}
and 
\begin{equation}t^2=1.\end{equation}
 
For any $i=1,\dots,n$, set
\begin{equation}\label{baxt-s}
\bs{i}(p,p',a,a'):=s_i+\frac{\delta_{p,p'}}{a-a'},
\end{equation}
where  $\delta_{p,p'}$ is the Kronecker delta.
For $p=p'$ the elements (\ref{baxt-s}) are called Baxterized elements; the parameters
$a$ and $a'$ are referred to as spectral parameters. We also define
\begin{equation}\label{baxt-t}
\bt(p):=\frac{1}{2}(1+p t).
\end{equation}
The following relation is satisfied and will be used later:
\begin{equation}\label{baxt-inv}
\bs{i}(p,p',a,a')\bs{i}(p',p,a',a)=\frac{(a-a')^2-\delta_{p,p'}}{(a-a')^2}\ \ \ \textrm{for $i=1,\dots,n$}.
\end{equation}

\vskip .2cm
Define the elements $j_i$, $i=1,\dots,n+1$, and $\tj_i$, $i=1,\dots,n+1$, of the group algebra $\mathbb{C}{\sf B_{n+1}}$ by the following initial conditions and recursions:
\begin{equation}\label{def-jm}
j_1=t\ ,\quad j_{i+1}=s_ij_is_i\quad\textrm{and}\quad \tj_1=0\ ,\quad \tj_{i+1}=s_i\tj_is_i+\frac{1}{2}(s_i+j_is_ij_i).
\end{equation}
The elements $j_i$ and $\tj_i$ are analogues, for the group ${\sf B}_{n+1}$, of the Jucys--Murphy elements.
The elements $j_i$, $i=1,\dots,n+1$, and $\tj_i$, $i=1,\dots,n+1$, form a maximal commutative set in $\mathbb{C}{\sf B}_{n+1}$, see \cite{OP,Pu,Ra}; in addition, $j_i$ and $\tj_i$ commute with all generators $s_k$, except $s_i$ and $s_{i-1}$:
\begin{equation}\label{com-jm}
j_is_k=s_kj_i\textrm{ and }\tj_is_k=s_k\tj_i\quad\textrm{if $k\neq i-1,i$.}
\end{equation}

\vskip .2cm
Let $\lambda\vdash n+1$ be a partition of length $n+1$, that is, $\lambda=(\lambda_1,\dots,\lambda_k)$, where $\lambda _j$, $j=1,\dots,k$, 
are positive integers, $\lambda_1\geqslant\lambda_2\geqslant\dots\geqslant\lambda_k$ and $n+1=\lambda_1+\dots+\lambda_k$.
We identify partitions with their Young diagrams: the Young diagram of $\lambda$ is a left-justified array of rows of
nodes containing $\lambda_j$ nodes in the $j$-th row, $j=1,\dots,k$; the rows are numbered from top to bottom.

\vskip .2cm
A $2$-partition, or a Young $2$-diagram, of length $n+1$ is a pair of partitions
such that the sum of their lengths 
equals $n+1$.
A $2$-node $\alpha^{(2)}$ is a pair $(\alpha,k)$ consisting of a usual node $\alpha$ and an integer $k=1,2$. The integer $k$ will be called {\it position}
of the $2$-node. A set of $2$-nodes can be equivalently 
described by an ordered pair of sets of nodes (the integer $k$ of a $2$-node $(\alpha,k)$ indicates to which set the node $\alpha$ belongs). A $2$-partition $\lambda^{(2)}$ is a set of $2$-nodes such that the subset consisting of the $2$-nodes having position $p$ is a usual partition, $p=1,2$.

\vskip .2cm
For a $2$-node $\alpha^{(2)}=(\alpha,k)$ lying in the line $x$ and the column $y$ of the 
$k$-th diagram, we denote by $\cc(\alpha^{(2)})$ the classical content of the node $\alpha$, $\cc(\alpha^{(2)}):=\cc(\alpha)=y-x$. Let $\{\xi_1,\xi_2\}$ be the  set of distinct square roots of unity,  ordered arbitrarily; we define  
also $p(\alpha^{(2)}):=\xi_k$.

\vskip .2cm
For a $2$-partition $\lambda^{(2)}$, a $2$-node $\alpha^{(2)}$ of $\lambda^{(2)}$ is called {\it removable} if the set of $2$-nodes obtained from $\lambda^{(2)}$ by removing $\alpha^{(2)}$ is still a $2$-partition. A $2$-node $\beta^{(2)}$ not in $\lambda^{(2)}$ is called {\it  addable} if the set of $2$-nodes obtained from $\lambda^{(2)}$ by adding $\beta^{(2)}$ is still a $2$-partition. For a $2$-partition $\lambda^{(2)}$, we denote by 
${\cal{E}}_-(\lambda^{(2)})$ the set of removable $2$-nodes of $\lambda^{(2)}$ and by ${\cal{E}}_+(\lambda^{(2)})$ the set of addable $2$-nodes of 
$\lambda^{(2)}$.

\vskip .2cm
Let $\lambda^{(2)}$ be a $2$-partition of length $n+1$. A standard $2$-tableau of shape $\lambda^{(2)}$ is obtained by placing the numbers $1,\dots,n+1$ in the $2$-nodes of the diagrams of $\lambda^{(2)}$ in such a way that the numbers in the nodes ascend along rows and columns in every diagram. For a standard $2$-tableau ${\cal{T}}$ of shape $\lambda^{(2)}$ let $\alpha^{(2)}_i$ be the $2$-node of ${\cal{T}}$ with number $i$, $i=1,\dots,n+1$; 
we set $\cc({\cal{T}}|i):=\cc(\alpha^{(2)}_i)$ and $p({\cal{T}}|i):=p(\alpha^{(2)}_i)$.

\vskip .2cm
The hook of a node $\alpha$ of a partition $\nu$ is the set of nodes of  $\nu$ consisting of the node $\alpha$ and the nodes which lie either under $\alpha$ in the same column or to the right of $\alpha$ in the same row; the hook length $h_{\nu}(\alpha)$ of $\alpha$ is the cardinality of 
the hook of $\alpha$. 
We extend this definition to $2$-nodes. For a $2$-node $\alpha^{(2)}=(\alpha,k)$ of a $2$-partition $\nu^{(2)}$, the hook length of $\alpha^{(2)}$ in $\nu^{(2)}$, which we denote by $h_{\nu^{(2)}}(\alpha^{(2)})$, is the hook length of the node $\alpha$ in the $k$-th partition of $\nu^{(2)}$.
Let
\begin{equation}\label{def-f}
f_{\nu^{(2)}}:=\bigl(\prod_{\alpha^{(2)}\in\nu^{(2)}}h_{\nu^{(2)}}(\alpha^{(2)})\bigr)^{-1}.
\end{equation}

\subsection{Idempotents of the group ${\sf B}_{n+1}$}\label{sub-idem}

The representation theory of the Coxeter groups of type $B$ was developed by A. Young \cite{Y1}, see also \cite{OP,Pu,Ra}. The irreducible representations of ${\sf B}_{n+1}$ are in bijection with the $2$-partitions of length $n+1$.
The elements of the semi-normal basis of the irreducible representation of ${\sf B}_{n+1}$ corresponding to the $2$-partition $\lambda^{(2)}$ are parameterized by the standard $2$-tableaux of shape $\lambda^{(2)}$. For a standard $2$-tableau ${\cal{T}}$, we denote by $E_{{\cal{T}}}$ the primitive idempotent of ${\sf B}_{n+1}$ corresponding to ${\cal{T}}$.
The Jucys--Murphy elements are diagonal in the semi-normal basis; moreover, we have, for any $i=1,\dots,n+1$,
\begin{equation}\label{spec-JM}
j_i E_{{\cal{T}}}=E_{{\cal{T}}} j_i=p_i E_{{\cal{T}}}\quad\textrm{and}\quad \tj_i E_{{\cal{T}}}=E_{{\cal{T}}} \tj_i=\cc_i E_{{\cal{T}}}.
\end{equation}
Here we set  $p_i:=p({\cal{T}}|i)$ and $\cc_i:=\cc({\cal{T}}|i)$ for all $i=1,\dots,n+1$ for brevity.
Due to the maximality of the commutative set $\{ j_i,\tj_i\}_{i=1,\dots,n+1}$ of Jucys--Murphy elements, the idempotent $E_{{\cal{T}}}$ can be expressed in terms 
of $j_i,\tj_i$, $i=1,\dots,n+1$.
Let $\alpha^{(2)}$ be the $2$-node of ${\cal{T}}$ with the number $n+1$. As the tableau ${\cal{T}}$ is standard, the $2$-node $\alpha^{(2)}$ of $\lambda^{(2)}$ is removable. Let ${\cal{U}}$ be the standard $2$-tableau obtained from ${\cal{T}}$ by removing the $2$-node $\alpha^{(2)}$ and let $\mu^{(2)}$ be the shape of ${\cal{U}}$.
The inductive formula for $E_{{\cal{T}}}$ in terms of the Jucys--Murphy elements reads:
\begin{equation}\label{idem-JM}
E_{{\cal{T}}}=E_{{\cal{U}}}\prod_{\beta^{(2)}\colon \begin{array}{l}\scriptstyle{\beta^{(2)}\in 
{\cal{E}}_+(\mu^{(2)})}\\\scriptstyle{\cc(\beta^{(2)})\neq \cc(\alpha^{(2)})}\end{array}  }\hspace{-0.6cm} \frac{\tj_{n+1}-\cc(\beta^{(2)})}{\cc(\alpha^{(2)})-\cc(\beta^{(2)})}\prod_{\beta^{(2)}\colon \begin{array}{l}\scriptstyle{\beta^{(2)}\in 
{\cal{E}}_+(\mu^{(2)})}\\\scriptstyle{p(\beta^{(2)})\neq p(\alpha^{(2)})}\end{array}  }\hspace{-0.5cm}\frac{j_{n+1}-p(\beta^{(2)})}{p(\alpha^{(2)})-p(\beta^{(2)})}\ .
\end{equation}
Note that the second product in the right hand side of (\ref{idem-JM}) contains only one term (it will not be so for the cyclotomic groups $G(m,1,n)$ with 
$m>2$).
We have ${\sf B}_0\cong \mathbb{C}$ and $E_{{\cal{U}}_0}=1$  for the (unique) $2$-tableau ${\cal{U}}_0$ of length $0$.

\vskip .2cm
Let $\{ {\cal{T}}_1,\dots,{\cal{T}}_k\}$ be the set of pairwise different standard $2$-tableaux that can be obtained from ${\cal{U}}$ by adding a $2$-node with the number $n+1$. The following formula:
\[E_{{\cal{U}}}=\sum_{i=1}^k E_{{\cal{T}}_i},\]
together with (\ref{spec-JM}) implies that the rational function
\[ E_{{\cal{U}}}\  \frac{u-\cc_{n+1}}{u-\tj_{n+1}}\ \frac{v-p_{n+1}}{v-j_{n+1}}\] 
is non-singular at $u=\cc_{n+1}$ and $v=p_{n+1}$, and, moreover, 
\begin{equation}\label{idem-JM-inter} E_{{\cal{U}}}\ 
\frac{u-\cc_{n+1}}{u-\tj_{n+1}}\ \frac{v-p_{n+1}}{v-j_{n+1}}\ 
\Bigr\rvert_{\hspace{-0.15cm}\textrm{\scriptsize{$\begin{array}{l}u=\cc_{n+1}\\[-0.15em]v=p_{n+1}\end{array}$}}}=E_{{\cal{T}}}\ .\end{equation}
Since $j_{n+1}$ takes values $\pm 1$, the  rational function $\frac{v-p_{n+1}}{v-j_{n+1}}$ is non-singular at $v=p_{n+1}$ and 
\begin{equation}\label{eq-v-j}
\frac{v-p_{n+1}}{v-j_{n+1}}\Bigr\rvert_{v=p_{n+1}}=\frac{1}{2}(1+p_{n+1}j_{n+1})\ .
\end{equation}
For the clarity of the calculations in the sequel, we define, generalizing (\ref{baxt-t}),
\begin{equation}\label{def-jbax}
\bj{i}(p ):=\frac{1}{2}(1+p j_i)\ \ \textrm{for $i=1,\dots,n+1$.}
\end{equation}
Combining (\ref{idem-JM-inter}) and (\ref{eq-v-j}), we obtain the following formula for the idempotent $E_{{\cal{T}}}$:
\begin{equation}\label{idem-JM2}
E_{{\cal{T}}}=E_{{\cal{U}}}\ \bj{n+1}(p_{n+1})\ \frac{u-\cc_{n+1}}{u-\tj_{n+1}}\Bigr\rvert_{u=\cc_{n+1}}.
\end{equation}

\section{{\hspace{-0.55cm}.\hspace{0.55cm}}Fusion formula for idempotents of ${\sf B}_{n+1}$}\label{sec-fus-B}

We start with the following Lemma which will be useful in the sequel.
\begin{lemm}
{\hspace{-.2cm}.\hspace{.2cm}}\label{lemm-tj}
For any integer $l$, $1\leqslant l\leqslant n$, we have
\begin{itemize}
 \item[{\rm (i)}] $\ \tj_{n+1}=s_ns_{n-1}\dots s_l\tj_ls_l\dots s_{n-1}s_n+\displaystyle{\frac{1}{2}\sum_{i=l}^n} s_n\dots s_{i+1}s_is_{i+1}\dots s_n(1+j_{n+1}j_i)$,
\item[{\rm (ii)}] $\ \bj{l}(p)s_l\dots s_{n-1}s_n\tj_{n+1}\!\!=\bj{l}(p)\tj_ls_l\dots s_{n-1}s_n+\displaystyle{\frac{1}{2}\sum_{i=l}^n} s_ls_{l+1}\dots s_{i-1}\cdot s_{i+1}\dots s_{n-1}s_n\,\bj{i}(p)(1+j_{n+1}j_i)$;
\end{itemize}
the product $s_ls_{l+1}\dots s_{i-1}$ in the right-hand-side of (ii) is understood to be equal to $1$ if $i= l$.
\end{lemm}
\emph{Proof.} We prove (i) by induction on $n-l$. The basis of the induction is, for $l=n$, the formula $\tj_{n+1}=s_n\tj_ns_n+\frac{1}{2}(s_n+s_nj_{n+1}j_n)$ which follows from the definition (\ref{def-jm}) of the Jucys--Murphy elements, namely from $\tj_{n+1}=s_n\tj_ns_n+\frac{1}{2}(s_n+j_ns_nj_n)$ and $j_ns_n=s_nj_{n+1}$. Now assume that
\[\tj_{n+1}=s_ns_{n-1}\dots s_{l+1}\tj_{l+1}s_{l+1}\dots s_{n-1}s_n+\frac{1}{2}\sum_{i=l+1}^n s_n\dots s_{i+1}s_is_{i+1}\dots s_n(1+j_{n+1}j_i),\]
and replace $\tj_{l+1}$ by $s_l\tj_ls_l+\frac{1}{2}(s_l+j_ls_lj_l)$. One obtains the assertion (i) using that $j_ls_{l+1}\dots s_{n-1}s_n=s_{l+1}\dots s_{n-1}s_nj_l$ and that $j_ls_ls_{l+1}\dots s_{n-1}s_n=s_ls_{l+1}\dots s_{n-1}s_nj_{n+1}$.

\vskip .2cm
Using (i), we find
\[\bj{l}(p)s_l\dots s_{n-1}s_n\tj_{n+1}=\bj{l}(p)\Bigl(\tj_ls_l\dots s_{n-1}s_n+\frac{1}{2}\sum_{i=l}^n s_ls_{l+1}\dots s_{i-1}\cdot 
s_{i+1}\dots s_{n-1}s_n(1+j_{n+1}j_i)\Bigr),\]
and (ii) follows since 
$j_ls_ls_{l+1}\dots s_{i-1}=s_ls_{l+1}\dots s_{i-1}j_i$ and $j_i$ commutes with $s_{i+1}\dots s_{n-1}s_n$.
\hfill$\square$

\vskip .2cm
Let $\phi_{1}(v,u):=\bt(v)$; for $k=1,\dots,n$ define
\begin{equation}\label{def-phi}\begin{array}{cc}\phi_{k+1}(v_1,\dots,v_k,v,u_1,\dots,u_k,u):=
\bs{k}(v,v_k,u,u_k)\phi_{k}(v_1,\dots,v_{k-1},v,u_1,\dots,u_{k-1},u)s_k&\\[.2cm]
=\bs{k}(v,v_k,u,u_k)\bs{k-1}(v,v_{k-1},u,u_{k-1})\dots \bs{1}(v,v_1,u,u_1)\bt(v)s_1\dots s_{k-1}s_k\ .&
\end{array}\end{equation}

\vskip .2cm
Define the following rational function with values in the group ring of ${\sf B}_{n+1}$:
\begin{equation}\label{def-Psi}
\Phi(v_1,\dots,v_{n+1},
u_1,\dots,u_{n+1}):=\prod_{k=0,\dots,n}^{\leftarrow}\phi_{k+1}(v_1,\dots,v_k,v_{k+1},u_1,\dots,u_k,u_{k+1})\ ;
\end{equation}
the arrow over $\prod$ indicates that the (non-commuting) factors are taken in the descending order.

\vskip .2cm
Let  $\lambda^{(2)}$ be a $2$-partition of length $n+1$ and ${\cal{T}}$ a standard $2$-tableau of shape $\lambda^{(2)}$. 
{}For $i=1,\dots,n+1$, set $p_i:=p({\cal{T}}|i)$ and $\cc_i:=\cc({\cal{T}}|i)$.

\begin{theo}
{\hspace{-.2cm}.\hspace{.2cm}}\label{prop-fus}
The idempotent $E_{{\cal{T}}}$ corresponding to the standard $2$-tableau ${\cal{T}}$ of shape $\lambda^{(2)}$ can be obtained by the following consecutive evaluations
\begin{equation}\label{eq-idem-fin}
E_{{\cal{T}}}=f_{\lambda^{(2)}}\Phi(v_1,\dots,v_{n+1},u_1,\dots,u_{n+1})\Bigr\rvert_{v_i=p_i,\ i=1,\dots,n+1}\Bigr\rvert_{u_1=\cc_1}\dots\Bigr\rvert_{u_n=\cc_n}\Bigr\rvert_{u_{n+1}=\cc_{n+1}}\ .
\end{equation}
\end{theo}
\emph{Proof.} 
Define
\begin{equation}\label{def-F}
F_{{\cal{T}}}(u):=\frac{u-\cc_{n+1}}{u}\prod_{i=1}^n\frac{(u-\cc_i)^2}{(u-\cc_i)^2-\delta_{p_i,p_{n+1}}}\ .
\end{equation} 
Let ${\cal{U}}$ be the standard $2$-tableau obtained from ${\cal{T}}$ by removing the $2$-node 
with the number $n+1$ and let $\mu^{(2)}$ be the shape of ${\cal{U}}$.

\begin{prop}{\hspace{-.2cm}.\hspace{.2cm}}\label{prop-idem1}
We have 
\begin{equation}\label{eq-idem1}
F_{{\cal{T}}}(u)\phi_{n+1}(p_1,\dots,p_n,p_{n+1},\cc_1,\dots,\cc_n,u)E_{{\cal{U}}}=\frac{u-\cc_{n+1}}{u-\tj_{n+1}}\ \bj{n+1}(p_{n+1})E_{{\cal{U}}}.
\end{equation}
\end{prop}
\emph{Proof.} We prove (\ref{eq-idem1}) by induction on $n$. As $\cc_1=0$ and $\tj_1=0$, the basis of induction for $n=0$ is the formula $\bt(p_1)=\bj{1}(p_1)$ which is satisfied by definition, see (\ref{baxt-t}) and (\ref{def-jbax}).

\vskip .2cm 
If $p_{n+1}\neq p_i$, $i=2,\dots,n$, then fix $l=1$. Otherwise fix $l$ such that $p_{n+1}=p_l$ and $p_{n+1}\neq p_i$, $i=l+1,\dots,n$.

Define ${\cal{V}}$ to be the standard $2$-tableau obtained from ${\cal{U}}$ by removing the $2$-nodes containing the numbers 
$l+1,\dots,n$ and ${\cal{W}}$ to be the standard $2$-tableau obtained from ${\cal{V}}$ by removing the $2$-node with the number $l$. We will use that $E_{{\cal{W}}}E_{{\cal{U}}}=E_{{\cal{U}}}$ and that $E_{{\cal{W}}}$ commutes with $s_i$, for $i=l,l+1,\dots,n$.
We rewrite the left-hand side of (\ref{eq-idem1}) as
\[F_{{\cal{T}}}(u)s_n\dots s_{l+1}\bs{l}(p_{n+1},p_l,u,\cc_l)\cdot\phi_l(p_1,\dots,p_{l-1},p_{n+1},\cc_1,\dots,\cc_{l-1},u) E_{{\cal{W}}}\cdot s_ls_{l+1}\dots s_n E_{{\cal{U}}}\ .\]
If $l>1$ then we use the induction hypothesis, namely \[\phi_l(p_1,\dots,p_{l-1},p_{l},\cc_1,\dots,\cc_{l-1},u) E_{{\cal{W}}}=(F_{{\cal{V}}}(u))^{-1}\frac{u-\cc_l}{u-\tj_l}\bj{l}(p_l)E_{{\cal{W}}}\ ,\]
and we notice that $p_{n+1}=p_l$.

\noindent If $l=1$ then $E_{{\cal{W}}}=1$, $F_{{\cal{V}}}(u)=1$ and, by definition, $\phi_1(p_{n+1},u)$ is equal to $\bj{1}(p_{n+1})$. Thus, in both situations (for $l=1$ and for $l>1$), we obtain for the left-hand-side of (\ref{eq-idem1}):
\[F_{{\cal{T}}}(u)(F_{{\cal{V}}}(u))^{-1}s_n\dots s_{l+1}\bs{l}(p_{n+1},p_l,u,\cc_l)\frac{u-\cc_l}{u-\tj_l}\bj{l}(p_{n+1})s_ls_{l+1}\dots s_n E_{{\cal{U}}}\ .\]
Therefore, the equality (\ref{eq-idem1}) is equivalent to
\begin{equation}\label{eq-inter1}
\begin{array}{l}
F_{{\cal{T}}}(u)(F_{{\cal{V}}}(u))^{-1}(u-\cc_l)\bj{l}(p_{n+1})s_ls_{l+1}\dots s_n(u-\tj_{n+1})E_{{\cal{U}}}\\[0.5em]
\hspace{1cm}=  \displaystyle{\frac{(u-\cc_{n+1})(u-\cc_l)^2}{(u-\cc_l)^2-\delta_{p_l,p_{n+1}}}}(u-\tj_l)\bs{l}(p_l,p_{n+1},\cc_l,u)s_{l+1}\dots s_n\bj{n+1}(p_{n+1})E_{{\cal{U}}},
\end{array}\end{equation}
where, in moving $s_n\dots s_{l+1}\bs{l}(p_{n+1},p_l,u,\cc_l)\displaystyle{\frac{1}{u-\tj_l}}$ to the right-hand-side and 
$\displaystyle{\frac{1}{u-\tj_{n+1}}}$ to the left-hand-side, we have used that $\tj_{n+1}$ commutes with $j_{n+1}$ and $E_{{\cal{U}}}$, and also the formula (\ref{baxt-inv}) to take the inverse of $\bs{l}(p_{n+1},p_l,u,\cc_l)$.

\vskip .2cm
To prove the equality (\ref{eq-inter1}), first notice that we have
\[F_{{\cal{T}}}(u)(F_{{\cal{V}}}(u))^{-1}(u-\cc_l)=(u-\cc_{n+1})\prod_{i=1}^n\frac{(u-\cc_i)^2}{(u-\cc_i)^2-\delta_{p_i,p_{n+1}}}\prod_{i=1}^{l-1}\left(\frac{(u-\cc_i)^2}{(u-\cc_i)^2-\delta_{p_i,p_{l}}}\right)^{-1},\]
which gives, since $p_i\neq p_{n+1}$ if $i>l$ and $p_l=p_{n+1}$ if $l>1$,
\begin{equation}\label{coeff-inter1}F_{{\cal{T}}}(u)(F_{{\cal{V}}}(u))^{-1}(u-\cc_l)=\frac{(u-\cc_{n+1})(u-\cc_l)^2}{(u-\cc_l)^2-\delta_{p_l,p_{n+1}}}.\end{equation}
So it remains to prove that
\begin{equation}\label{eq-inter2}
\bj{l}(p_{n+1})s_ls_{l+1}\dots s_n(u-\tj_{n+1})E_{{\cal{U}}}=(u-\tj_l)\bs{l}(p_l,p_{n+1},\cc_l,u)s_{l+1}\dots s_n\bj{n+1}(p_{n+1})E_{{\cal{U}}}.
\end{equation}
Expand $\bs{l}(p_l,p_{n+1},\cc_l,u)$
in the right hand side of (\ref{eq-inter2}): 
\begin{equation}\label{eq-inter2'}\bigl((u-\tj_l)s_l-\delta_{p_l,p_{n+1}}\frac{\tj_l-u}{\cc_l-u}\bigr)s_{l+1}\dots s_n\bj{n+1}(p_{n+1})E_{{\cal{U}}}\ .\end{equation}
As $\tj_l$ commutes with $s_{l+1}\dots s_n$ and $j_{n+1}$ and $\tj_lE_{{\cal{U}}}=\cc_lE_{{\cal{U}}}$, we find that the expression (\ref{eq-inter2'}) equals
\[\Bigl((u-\tj_l)s_ls_{l+1}\dots s_n\bj{n+1}(p_{n+1})-\delta_{p_l,p_{n+1}}s_{l+1}\dots s_n\bj{n+1}(p_{n+1})\Bigr)E_{{\cal{U}}}\ .\]
Then using that $s_ls_{l+1}\dots s_nj_{n+1}=j_ls_ls_{l+1}\dots s_n$, we obtain for the right-hand-side of (\ref{eq-inter2}):
\begin{equation}\label{eq-inter2-d}
\bigl((u-\tj_l)\bj{l}(p_{n+1})s_ls_{l+1}\dots s_n-\delta_{p_l,p_{n+1}}s_{l+1}\dots s_n\bj{n+1}(p_{n+1})\bigr)E_{{\cal{U}}}.
\end{equation}
Using the Lemma \ref{lemm-tj}, (ii), we write the left hand side of (\ref{eq-inter2}) in the form
\begin{equation}\label{eq-inter3}
\Bigl(\bj{l}(p_{n+1})(u-\tj_{l})s_ls_{l+1}\dots s_n-\frac{1}{2}\sum_{i=l}^ns_ls_{l+1}\dots s_{i-1}\cdot s_{i+1}\dots s_{n-1}s_n\bj{i}(p_{n+1})(1+j_{n+1}j_i)\Bigr)E_{{\cal{U}}}\ .
\end{equation}
As $j_iE_{{\cal{U}}}=p_iE_{{\cal{U}}}$, $i=1,\dots,n$, the expression (\ref{eq-inter3}) is equal to
\begin{equation}\label{eq-inter4}\Bigl(\bj{l}(p_{n+1})(u-\tj_{l})s_ls_{l+1}\dots s_n-\frac{1}{2}\sum_{i=l}^ns_ls_{l+1}\dots s_{i-1}\cdot s_{i+1}\dots s_{n-1}s_n\frac{1}{2}(1+p_ip_{n+1})(1+p_ij_{n+1})\Bigr)E_{{\cal{U}}}\ .\end{equation}
Since $p_kp_{n+1}=-1$, $k=l+1,\dots,n$, we finally obtain for the left-hand-side of (\ref{eq-inter2}):
\begin{equation}\label{eq-inter2-g}
\Bigl(\bj{l}(p_{n+1})(u-\tj_{l})s_ls_{l+1}\dots s_n-\delta_{p_l,p_{n+1}}s_{l+1}\dots s_{n-1}s_n\bj{n+1}(p_l)\bigr)E_{{\cal{U}}}\ .
\end{equation}
The comparison of (\ref{eq-inter2-g}) and (\ref{eq-inter2-d}) proves the equality (\ref{eq-inter2}). 
\hfill $\square$

\begin{prop}
{\hspace{-.2cm}.\hspace{.2cm}}\label{prop-f}
The rational function $F_{{\cal{T}}}(u)$, defined by (\ref{def-F}), is regular at $u=\cc_{n+1}$ and moreover
\begin{equation}\label{eq-F}
F_{{\cal{T}}}(\cc_{n+1})=f_{\lambda^{(2)}}(f_{\mu^{(2)}})^{-1}.
\end{equation}
\end{prop}
\emph{Proof.} The Proposition will directly follow from the result, used in \cite{Mo}, concerning usual tableaux. Let $\lambda$ be a partition of length $m+1$, with $m\leq n$. let $S$ be a subset of $\{1,\dots,n+1\}$ such that $S$ contains the number $n+1$ and has a cardinal equal to $m+1$. Let $\tilde{{\cal{T}}}$ be a tableau of shape $\lambda$ filled with numbers belonging to $S$ such that the numbers in the nodes are in strictly ascending orders along rows and columns in right and down directions. Let $\gamma$ be the node of $\tilde{{\cal{T}}}$ with the number $n+1$ and $\mu$ be the 
shape of the tableau obtained from $\tilde{{\cal{T}}}$ by removing the node $\gamma$.
Define the following rational function
\begin{equation}\label{def-tF}\tilde{F}_{\tilde{{\cal{T}}}}(u)=\frac{u-\cc(\gamma)}{u}\prod_{\alpha\in\mu}\frac{\left(u-\cc(\alpha)\right)^2}{(u-\cc(\alpha)+1)(u-\cc(\alpha)-1)}.\end{equation}
The product in the right hand side of (\ref{def-tF}) depends only on the 
shape $\mu$ and one has
\[ \prod_{\alpha\in\mu}\frac{\left(u-\cc(\alpha)\right)^2}{(u-\cc(\alpha)+1)(u-\cc(\alpha)-1)}=u\ \displaystyle{\prod_{\beta\in {\cal{E}}_-(\mu)}\left(u-\cc(\beta)\right)}\displaystyle{\prod_{\alpha\in {\cal{E}}_+(\mu)}\left(u-\cc(\alpha)\right)^{-1}}\ ,
\]
where ${\cal{E}}_-(\mu)$ (respectively, ${\cal{E}}_+(\mu)$) is the set of removable (respectively, addable) nodes of $\mu$.
Therefore the rational function $\tilde{F}_{\tilde{{\cal{T}}}}(u)$ is regular at $u=\cc(\gamma)$ and moreover
\begin{equation}\label{tF} \tilde{F}_{\tilde{{\cal{T}}}}(\cc(\gamma))=\displaystyle{\prod_{\beta\in {\cal{E}}_-(\mu)}\left(\cc(\gamma)-\cc(\beta)\right)}\displaystyle{\prod_{\alpha\in {\cal{E}}_+(\mu)\backslash\{\gamma\}}\left(\cc(\gamma)-\cc(\alpha)\right)^{-1}}\ .
\end{equation}
It is known that the right hand side of (\ref{tF}) is equal to 
\begin{equation}\label{eq-lemm-f}\prod_{\alpha\in\lambda}\bigl(h_{\lambda}(\alpha)\bigr)^{-1}\prod_{\alpha\in\mu}h_{\mu}(\alpha).\end{equation}

\vskip .2cm
Define $\tilde{{\cal{T}}}$ to be the tableau of the standard $2$-tableau ${\cal{T}}$ which contains the node with number $n+1$. 
The assertion of the Proposition \ref{prop-f} follows, the only observation one has to make is that the $2$-nodes $(\alpha,k)$ with $p_k\neq p_{n+1}$ do not contribute to (\ref{eq-F}). 
\hfill $\square$

\vskip .2cm
The Theorem \ref{prop-fus} follows, by induction on $n$, from the formula (\ref{idem-JM2}), the Proposition \ref{prop-idem1} and the Proposition \ref{prop-f}.\hfill $\square$

\section{{\hspace{-0.55cm}.\hspace{0.55cm}}Fusion procedure for the complex reflection group $G(m,1,n+1)$}\label{sec-fus-G} 

We extend the results of the previous Section to the complex reflection groups $G(m,1,n+1)$ for all positive integers $m$. We skip the proofs when they are completely similar to the proofs in the preceding Section; we only indicate modifications.

\subsection{Definitions}

The complex reflection group $G(m,1,n+1)$ is generated by the elements $s_1,\dots,s_n$ and $t$ with the defining relations (\ref{rel-def-An}), 
(\ref{rel-def-Bn}) and
\begin{equation}\label{rel-def-Gm}
t^m=1. 
\end{equation} 
In particular, $G(1,1,n+1)$ is isomorphic to the symmetric group $S_{n+1}$ and $G(2,1,n+1)$ to ${\sf B}_{n+1}$.

\vskip .2cm
We extend the definition (\ref{baxt-s}) to the generators $s_1,\dots,s_n$ of $G(m,1,n+1)$:
\begin{equation}\label{baxt-s-m}
\bs{i}(p,p',a,a'):=s_i+\frac{\delta_{p,p'}}{a-a'},\ \ i=1,\dots,n\ .
\end{equation}

The Jucys--Murphy elements for the group $G(m,1,n+1)$ are the elements $j_i$, $i=1,\dots,n+1$, and $\tj_i$, $i=1,\dots,n+1$, of the group ring defined inductively by the following initial conditions and recursions:
\begin{equation}\label{def-jm-m}
j_1=t\ ,\quad j_{i+1}=s_ij_is_i\quad\text{and}\quad \tj_1=0\ ,\quad \tj_{i+1}=s_i\tj_is_i+\frac{1}{m}\sum\limits_{k=0}^{m-1}j_i^ks_ij_i^{m-k}.
\end{equation}

For $m=1$, that is, for $S_{n+1}$, $j_k=1$, $k=1,\dots,n+1$; the recursion formula for $\tj_{i+1}$ reduces to $\tj_{i+1}=s_i\tj_i s_i+s_i$.

\vskip .3cm
As for
$m=2$, the elements $j_i$, $i=1,\dots,n+1$, and $\tj_i$, $i=1,\dots,n+1$, form a maximal commutative set 
in $\mathbb{C}G(m,1,n+1)$,
see \cite{OP,Pu}; in addition, $j_i$ and $\tj_i$ commute with all
$s_k$, except $s_i$ and $s_{i-1}$.

\vskip .2cm
The definitions of a $2$-partition, $2$-node, standard $2$-tableaux and hook length generalize naturally
to any $m$; for example, an $m$-partition is an $m$-tuple of partitions and an $m$-node $\alpha^{(m)}$ is a pair $(\alpha,k)$, $k=1,\dots ,m$.
For an $m$-node $\alpha^{(m)}=(\alpha,k)$ of an $m$-partition $\lambda^{(m)}$ such that the node $\alpha$ lies in the line $x$ and the column $y$ of the $k$-th diagram, we define $\cc(\alpha^{(m)}):=\cc(\alpha)=y-x$. Let $\{\xi_1,\dots,\xi_m\}$ be the  set of distinct $m$-th roots of unity, ordered arbitrarily; we define  
$p(\alpha^{(m)}):=\xi_k$.

\vskip .2cm 
As for $m=2$, we define
\begin{equation}\label{def-f-m}
f_{\nu^{(m)}}:=\bigl(\prod_{\alpha^{(m)}\in\nu^{(m)}}h_{\nu^{(m)}}(\alpha^{(m)})\bigr)^{-1}\ ,
\end{equation}
where $h_{\nu^{(m)}}(\alpha^{(m)})$ is the hook length of $\alpha^{(m)}$ calculated in $\nu^{(m)}$.

\vskip .2cm
The irreducible representations of $G(m,1,n+1)$ are parameterized by the $m$-partitions of length $n+1$; for a given $m$-partition $\lambda^{(m)}$, elements of the semi-normal basis of the corresponding representation are indexed by the standard $m$-tableaux of shape $\lambda^{(m)}$. We denote by $E_{{\cal{T}}}$ the idempotent of the group ring corresponding to a standard $m$-tableau ${\cal{T}}$.

\subsection{Fusion formula for idempotents of $G(m,1,n+1)$}

Let $\lambda^{(m)}$ be an $m$-partition of length $n+1$; fix a standard $m$-tableau ${\cal{T}}$ of shape $\lambda^{(m)}$. Let ${\cal{U}}$ be the standard $m$-tableau obtained by removing from ${\cal{T}}$ the $m$-node containing $n+1$ and let $\mu^{(m)}$ be the shape of ${\cal{U}}$. Set, for $i=1,\dots,n+1$, $p_i:=p({\cal{T}}|i)$ and $\cc_i:=\cc({\cal{T}}|i)$ 

\vskip .2cm
The same reasoning as in Subsection \ref{sub-idem} leads to the formula for $E_{{\cal{T}}}$ (cf  (\ref{idem-JM-inter})):
\begin{equation}\label{idem-JM-inter-m} E_{{\cal{U}}}\ 
\frac{u-\cc_{n+1}}{u-\tj_{n+1}}\ \frac{v-p_{n+1}}{v-j_{n+1}}\ 
\Bigr\rvert_{\hspace{-0.15cm}\textrm{\scriptsize{$\begin{array}{l}u=\cc_{n+1}\\[-0.15em]v=p_{n+1}\end{array}$}}}=E_{{\cal{T}}}\ .\end{equation}
Since $j_{n+1}$ takes values in $\{\xi_1,\dots,\xi_m\}$, the rational function $\frac{v-p_{n+1}}{v-j_{n+1}}$ is non-singular for $v=p_{n+1}$ and 
\begin{equation}\label{eq-v-j-m}
\frac{v-p_{n+1}}{v-j_{n+1}}\Bigr\rvert_{v=p_{n+1}}=\frac{1}{m}\sum\limits_{k=0}^{m-1}p_{n+1}^{m-k}j_{n+1}^k\ .
\end{equation}
An analogue of (\ref{def-jbax}) is
\begin{equation}\label{def-jbax-m}
\bj{i}(p):=\frac{1}{m}\sum\limits_{k=0}^{m-1}p^{m-k}j_i^k\ \ \textrm{for $i=1,\dots,n+1$.}
\end{equation}
For $i=1$ we shall write $\bt(p):=\frac{1}{m}\sum\limits_{k=0}^{m-1}p^{m-k}t^k$ instead of $\bj{1}(p)$.
Combining (\ref{idem-JM-inter-m}) and (\ref{eq-v-j-m}), we obtain  for the idempotent $E_{{\cal{T}}}$ (cf (\ref{idem-JM2})) 
\begin{equation}\label{idem-JM2-m}
E_{{\cal{T}}}=E_{{\cal{U}}}\ \bj{n+1}(p_{n+1})\ \frac{u-\cc_{n+1}}{u-\tj_{n+1}}\Bigr\rvert_{u=\cc_{n+1}}.
\end{equation}

\vskip .2cm
We generalize the Lemma \ref{lemm-tj} to an arbitrary positive integer $m$.
\begin{lemm}
{\hspace{-.2cm}.\hspace{.2cm}}\label{lemm-tj-m}
For any integer $l$, $1\leqslant l\leqslant n$, we have

\noindent
{\rm (i)} $\tj_{n+1}=s_ns_{n-1}\dots s_l\tj_ls_l\dots s_{n-1}s_n+\displaystyle{\frac{1}{m}}\sum\limits_{i=l}^n s_n\dots s_{i+1}s_is_{i+1}\dots s_n\sum\limits_{k=0}^{m-1}j_{n+1}^kj_i^{m-k}$.

\noindent
{\rm (ii)} $\bj{l}(p)s_l\dots s_{n-1}s_n\tj_{n+1}\!\!=\!\bj{l}(p)\tj_ls_l\dots s_{n-1}s_n\!+\!\displaystyle{\frac{1}{m}}\!\sum\limits_{i=l}^n s_ls_{l+1}\dots s_{i-1}\cdot
s_{i+1}\dots s_{n-1}s_n\,\bj{i}(p)\!\sum\limits_{k=0}^{m-1}j_{n+1}^kj_i^{m-k}$;

\noindent
the product $s_ls_{l+1}\dots s_{i-1}$ in the right-hand-side of (ii) is understood to be equal to $1$ if $i= l$.
\end{lemm}
\emph{Proof.} The proof is completely similar to the proof of the Lemma \ref{lemm-tj}.\hfill$\square$

\vskip .2cm
Let $\phi_{1}(v,u):=\bt(v)$; for $k=1,\dots,n$ define
\begin{equation}\label{def-phi-m}\begin{array}{cc}\phi_{k+1}(v_1,\dots,v_k,v,u_1,\dots,u_k,u):=
\bs{k}(v,v_k,u,u_k)\phi_{k}(v_1,\dots,v_{k-1},v,u_1,\dots,u_{k-1},u)s_k&\\[.2cm]
=\bs{k}(v,v_k,u,u_k)\bs{k-1}(v,v_{k-1},u,u_{k-1})\dots \bs{1}(v,v_1,u,u_1)\bt(v)s_1\dots s_{k-1}s_k\ .&
\end{array}\end{equation}
The formula (\ref{def-phi-m}) reads in the same way for any $m$; only the definition of $\bt(v)$ depends on $m$.

\vskip .2cm
Define the following rational function with values in the group ring of $G(m,1,n)$:
\begin{equation}\label{def-Psi-m}
\Phi(v_1,\dots,v_{n+1},u_1,\dots,u_{n+1}):=\prod_{k=0,\dots,n}^{\leftarrow}\phi_{k+1}(v_1,\dots,v_k,v_{k+1},u_1,\dots,u_k,u_{k+1})\ ;
\end{equation}
the arrow over $\prod$ indicates that the (non-commuting) factors are taken in the descending order.

\begin{theo}
{\hspace{-.2cm}.\hspace{.2cm}}\label{prop-fus-m}
The idempotent $E_{{\cal{T}}}$ corresponding to the standard $m$-tableau ${\cal{T}}$ of shape $\lambda^{(m)}$ can be obtained by the following consecutive evaluations
\begin{equation}\label{eq-idem-fin-m}
E_{{\cal{T}}}=f_{\lambda^{(m)}}\Phi(v_1,\dots,v_{n+1},u_1,\dots,u_{n+1})\Bigr\rvert_{v_i=p_i,\ i=1,\dots,n+1}\Bigr\rvert_{u_1=\cc_1}\dots\Bigr\rvert_{u_n=\cc_n}\Bigr\rvert_{u_{n+1}=\cc_{n+1}}\ .
\end{equation}
\end{theo}
\emph{Proof.} Define
\begin{equation}\label{def-F-m}
F_{{\cal{T}}}(u):=\frac{u-\cc_{n+1}}{u}\prod_{i=1}^n\frac{(u-\cc_i)^2}{(u-\cc_i)^2-\delta_{p_i,p_{n+1}}}\ .
\end{equation}
\begin{prop}
{\hspace{-.2cm}.\hspace{.2cm}}\label{prop-idem1-m}
We have 
\begin{equation}\label{eq-idem1-m}
F_{{\cal{T}}}(u)\phi_{n+1}(p_1,\dots,p_n,p_{n+1},\cc_1,\dots,\cc_n,u)E_{{\cal{U}}}=\frac{u-\cc_{n+1}}{u-\tj_{n+1}}\ \bj{n+1}(p_{n+1})E_{{\cal{U}}}.
\end{equation}
\end{prop}
\emph{Proof.} The proof follows the same lines as the proof of the Proposition \ref{prop-idem1}; actually it is exactly the same until the calculation of $\bj{l}(p_{n+1})s_ls_{l+1}\dots s_n(u-\tj_{n+1})E_{{\cal{U}}}$ just after the formula (\ref{eq-inter2-d}). We give the modified end of the proof.

Here, we rewrite $\bj{l}(p_{n+1})s_ls_{l+1}\dots s_n(u-\tj_{n+1})E_{{\cal{U}}}$ using the Lemma \ref{lemm-tj-m}, (ii):
\begin{equation}\label{eq-inter3-m}
\Bigl(\bj{l}(p_{n+1})(u-\tj_{l})s_ls_{l+1}\dots s_n-\frac{1}{m}\sum_{i=l}^ns_ls_{l+1}\dots s_{i-1}\cdot s_{i+1}\dots s_{n-1}s_n\bj{i}(p_{n+1})\sum\limits_{k=0}^{m-1}j_{n+1}^kj_i^{m-k}\Bigr)E_{{\cal{U}}}\ .
\end{equation}
As $j_iE_{{\cal{U}}}=p_iE_{{\cal{U}}}$ for $i=1,\dots,n$, the expression (\ref{eq-inter3-m}) is equal to
\begin{equation}\label{eq-inter4-m}\Bigl(\bj{l}(p_{n+1})(u-\tj_{l})s_ls_{l+1}\dots s_n-\frac{1}{m}\sum_{i=l}^ns_ls_{l+1}\dots s_{i-1}\cdot s_{i+1}\dots s_{n-1}s_n\sum\limits_{k=0}^{m-1}p_{n+1}^kp_i^{m-k}\bj{n+1}(p_i)\Bigr)E_{{\cal{U}}}\ .\end{equation}
Since 
$\frac{1}{m}\sum\limits_{k=0}^{m-1}p_{n+1}^kp_i^{m-k}=
\delta_{p_i,p_{n+1}}$
we obtain that $\bj{l}(p_{n+1})s_ls_{l+1}\dots s_n(u-\tj_{n+1})E_{{\cal{U}}}$ equals
\begin{equation}\label{eq-inter2-g-m}
\Bigl(\bj{l}(p_{n+1})(u-\tj_{l})s_ls_{l+1}\dots s_n-\delta_{p_l,p_{n+1}}s_{l+1}\dots s_{n-1}s_n\bj{n+1}(p_l)\bigr)E_{{\cal{U}}}\ .
\end{equation}
This concludes the proof.
\hfill $\square$

\vskip .2cm
The analogue of the Proposition \ref{prop-f} holds as well.
\begin{prop}
{\hspace{-.2cm}.\hspace{.2cm}}\label{prop-f-m}
The rational function $F_{{\cal{T}}}(u)$, defined by (\ref{def-F-m}), is regular at $u=\cc_{n+1}$ and moreover
\begin{equation}\label{eq-F-m}
F_{{\cal{T}}}(\cc_{n+1})=f_{\lambda^{(m)}}(f_{\mu^{(m)}})^{-1}. 
\end{equation}
\end{prop}
\emph{Proof.} The proof is completely similar to the proof of the Proposition \ref{prop-f}.\hfill $\square$

\vskip .2cm
Similarly to the Theorem \ref{prop-fus}, the Theorem \ref{prop-fus-m} follows, using induction on $n$, from the formula (\ref{idem-JM2-m}), the Proposition \ref{prop-idem1-m} and the Proposition \ref{prop-f-m}.\hfill $\square$

\vskip .2cm
{}For calculations, it is sometimes useful to write the function $\Phi$ in a slightly different form. Namely, let $\tilde{\phi}_{1}(v,u):=1$ and define
\begin{equation}\label{def-phi-m2}\begin{array}{cc}\tilde{\phi}_{k+1}(v_1,\dots,v_k,v,u_1,\dots,u_k,u):=
\bs{k}(v,v_k,u,u_k)\tilde{\phi}_{k}(v_1,\dots,v_{k-1},v,u_1,\dots,u_{k-1},u)s_k&\\[.2cm]
=\bs{k}(v,v_k,u,u_k)\bs{k-1}(v,v_{k-1},u,u_{k-1})\dots \bs{1}(v,v_1,u,u_1)s_1\dots s_{k-1}s_k\ ,&
\end{array}\end{equation}
for $k=1,\dots,n$. The elements $\tilde{\phi}_{k+1}(v_1,\dots,v_k,v,u_1,\dots,u_k,u)$ do not involve the generator $t$ and 
$\Phi(v_1,\dots,v_{n+1},u_1,\dots,u_{n+1})$, defined in (\ref{def-Psi-m}), equals 
\begin{equation}\label{def-Psi-m2}
\prod_{k=0,\dots,n}^{\leftarrow}\tilde{\phi}_{k+1}(v_1,\dots,v_{k+1},u_1,\dots,u_{k+1})
\cdot\bj{1}(v_1)\bj{2}(v_2)\dots\bj{n+1}(v_{n+1})\ .
\end{equation}
For example, let $m=2$; choose the order $\{ 1,-1\}$ on the set of square roots of 1. The primitive idempotent, corresponding to the 
standard 2-tableau ${\textrm{\tiny{$\left(
\fbox{\scriptsize{$1$}}\fbox{\scriptsize{$3$}}\, ,\fbox{\scriptsize{$2$}}\right)$}}}$ reads $s_2(1+s_1)s_2\bj{1}(1)\bj{2}(-1)\bj{3}(1)/16$.

\end{document}